\documentclass[twocolumn]{article}
\usepackage{pdfsync}
\usepackage{fullpage}
\usepackage[latin1]{inputenc}
\usepackage{color}
\usepackage{subfigure}
\usepackage{cite}
\usepackage{amsmath,amsfonts,amssymb,amsthm}

\DeclareMathOperator{\QFT}{QFT}
\usepackage{url,xspace}
\usepackage{graphicx}
\newtheorem{definition}{Definition}

\newtheorem{property}{Property}
\renewcommand{\v}[1]{\boldsymbol{#1}} % Vectors
\renewcommand{\i}{\ensuremath{\v{i}}\xspace}
\renewcommand{\j}{\ensuremath{\v{j}}\xspace}
\renewcommand{\k}{\ensuremath{\v{k}}\xspace}

\newcommand{\norm}[1]{\ensuremath{||#1||}\xspace}

\newcommand{\R}{\ensuremath{{\mathbb R}}}
\newcommand{\C}{\ensuremath{{\mathbb C}}}
\renewcommand{\H}{\ensuremath{{\mathbb H}}}
\newcommand*{\qconjugate}[1]{\overline{#1}}                % Quaternion conjugate.
\newcommand*{\cconjugate}[1]{{#1^\star}}                   % Complex    conjugate.
\DeclareMathOperator{\sign}{sign}
\newcommand{\dt}{\mathrm{d}t}

\newcommand{\dtau}{\mathrm{d}\tau}

\newcommand*{\hilbert}[1]{\ensuremath{\mathcal{H}\left[#1\right]}\xspace}
\newcommand*{\qhilbert}[2]{\ensuremath{\mathcal{H}_{#1}\left[#2\right]}\xspace}

\newcommand*{\iqfourier}[2]{\ensuremath{\mathcal{F}_{#1}^{-1}\left[#2\right]}\xspace}
\let\oldmu\mu\renewcommand{\mu}{\v{\oldmu}} % These macros should give a bold version
\let\oldxi\xi\renewcommand{\xi}{\v{\oldxi}}

\begin{document}

\title{Instantaneous frequency and amplitude of complex signals based on quaternion Fourier transform}

\author{Nicolas Le~Bihan\thanks{Nicolas Le~Bihan is with the CNRS,
                                GIPSA-Lab, D\'epartement Images et Signal,
                                961 Rue de la Houille Blanche,
                                Domaine Universitaire, BP~46,
                                38402~Saint Martin d'H\`eres cedex,
                                France,
                                email: nicolas.le-bihan@gipsa-lab.grenoble-inp.fr}
\and Stephen J. Sangwine\thanks{Stephen J. Sangwine is with the School of Computer
                                Science and Electronic Engineering, University of Essex,
                                Colchester, CO4 3SQ, United Kingdom,
                                email: sjs@essex.ac.uk}
\and Todd A. Ell\thanks{Todd A. Ell is with Goodrich Sensors \& Integrated Systems,
                        14300 Judicial Road, Burnsville, MN 55306 USA,
                        email: t.ell@ieee.org}
}

\maketitle

\begin{abstract}
The ideas of instantaneous amplitude and phase are well understood
for signals with real-valued samples, based on the analytic signal
which is a complex signal with one-sided Fourier transform.
We extend these ideas to signals with complex-valued samples,
using a quaternion-valued equivalent of the analytic signal
obtained from a one-sided quaternion Fourier transform which
we refer to as the \emph{hypercomplex representation} of the complex signal.
We present the necessary properties of the quaternion Fourier
transform, particularly its symmetries in the frequency domain
and formulae for convolution and the quaternion Fourier transform
of the Hilbert transform.
The hypercomplex representation may be interpreted as an ordered pair
of complex signals or as a quaternion signal.
We discuss its derivation and properties and show that its
quaternion Fourier transform is one-sided.
It is shown how to derive from the hypercomplex representation a
complex envelope and a phase.

A classical result in the case of real signals is that an
amplitude modulated signal may be analysed into its envelope
and carrier using the analytic signal provided that the modulating
signal has frequency content not overlapping with that of the
carrier.
We show that this idea extends to the complex case,
provided that the complex signal modulates an \emph{orthonormal}
complex exponential.
Orthonormal complex modulation can be represented mathematically
by a polar representation of quaternions previously derived by
the authors.
As in the classical case, there is a restriction of non-overlapping
frequency content between the modulating complex signal and the
orthonormal complex exponential.
We show that, under these conditions, modulation in the time
domain is equivalent to a frequency shift in the quaternion
Fourier domain.
Examples are presented to demonstrate these concepts.
\end{abstract}

\section{Introduction}
\label{Introsection}

The instantaneous amplitude and phase of a real-valued signal have been
understood since 1948 from the work of Ville \cite{Ville:1948}
and Gabor \cite{Gabor:1946} based on the \emph{analytic signal}.
A critical discussion of instantaneous amplitude and phase was given by
Picinbono in 1997,
particularly with reference to amplitude and phase modulation \cite[§V]{Picinbono:1997}.
The analytic signal and its associated modulation concepts can be simply described,
even though a full theoretical treatment is quite deep.
Given a real-valued signal $f(t)$,
the corresponding analytic signal $a(t)$ is a complex signal with real part equal to $f(t)$
and imaginary part orthogonal to $f(t)$.
The imaginary part is sometimes known as the quadrature signal
-- in the case where $f(t)$ is a sinusoid,
the imaginary part of the analytic signal is in quadrature,
that is with a phase difference of $-\pi/2$.
The orthogonal signal is related to $f(t)$ by the Hilbert transform
\cite{Hahn:1996,Hahn:Poularikas}.
The analytic signal has the interesting property that its modulus $|a(t)|$ is an
envelope of the signal $f(t)$.
The envelope is also known as the \emph{instantaneous amplitude}.
Thus if $f(t)$ is an amplitude-modulated sinusoid,
the envelope $|a(t)|$,
\emph{subject to certain conditions on the frequency content},
is the modulating signal.
The argument of the analytic signal,
$\angle\,a(t)$ is known as the \emph{instantaneous phase}
and its derivative is known as the \emph{instantaneous frequency}.
The analytic signal has a further very interesting property:
it has a one-sided Fourier transform.
Thus a simple method for constructing the analytic signal
(algebraically or numerically)
is to compute the Fourier transform of $f(t)$,
multiply the Fourier transform by a unit step which is zero for negative frequencies,
and then construct the inverse Fourier transform of the resulting spectrum.

In this paper we extend the above ideas to signals with complex-valued samples.
We show that if a quaternion Fourier transform is computed from a complex-valued
signal, and negative frequencies in the frequency domain representation are suppressed,
the inverse quaternion Fourier transform gives an equivalent of the analytic signal,
which we refer to in this paper as a \emph{hypercomplex representation} of the
complex signal.
Just as the classical analytic signal represents a real-valued signal by a complex
signal (that is with a pair of real signals),
the hypercomplex representation presented here represents a complex-valued signal
by a quaternion signal (that is with a pair of complex signals, based on the
Cayley-Dickson form of a quaternion).

We consider signals with complex-valued samples in a most general sense,
without any restrictions on the relationship between the real and imaginary parts.
However, we exclude the special case of an \emph{analytic signal}
(where the real and imaginary parts are orthogonal)
because it does not possess an interesting hypercomplex representation:
the two additional quaternion components of the hypercomplex representation
are simply copies of the two original components.
The case where the real and imaginary parts are correlated in some statistical
sense is called \emph{improper}, whereas an analytic signal is \emph{proper}.
It has been customary to handle such signals by using an augmented representation
consisting of a vector containing the signal and its conjugate \cite{Schreier:2010}.
Recently, Lilly and Olhede \cite{10.1109/TSP.2009.2031729} have published a paper
on bivariate analytic signal concepts.
Their approach is linked to a specific signal model,
the \emph{modulated elliptical signal},
which they illustrate with the example of a drifting oceanographic float,
and they utilise pairs of analytic signals derived from the pair of real
signals representing the time-varying positional coordinates of a drifting float.
In this paper a different approach, using hypercomplex representation,
is considered because it is more suited to the analysis of modulation,
and also because we have at our disposal the existing tools of quaternion Fourier
transforms, which provide a natural fit with modulation of a complex signal
by a complex signal, for reasons which should be evident later in the paper. Note that, in essence, the approach proposed in this paper is inspired by the work of B\"ulow \cite{Bulow:1999,BulowSommer:2001} and Labunets \cite{Labunets:2001} who considered Clifford Fourier transforms (quaternion in the case of B\"ulow) to process multidimensional signals. While B\"ulow considered 2D signals (images), here we consider 1D signals with 2D samples (complex valued samples).
We are not focusing on any particular application in this paper,
rather we propose a new approach to the processing of complex valued signals
(including proper and improper signals,
even though proper signals are well-described using the classical complex Fourier transform
--- they are a special case.)
Schreier and Scharf's book gives examples of applications \cite{Schreier:2010}.

The idea of amplitude modulation is extended to the case where the modulating
signal is complex,
and the `carrier' is an orthonormal complex exponential,
that is a complex exponential which is in a complex plane orthogonal to the
complex plane of the modulating signal
-- this concept is easily realised using quaternions,
since quaternion algebra has a 4-dimensional basis with
three mutually orthogonal imaginary units.
Subject to the same restrictions on frequency content as in the real-valued case,
amplitude modulation in the time domain corresponds to a frequency shift in the Fourier
domain.
Further, a complex envelope of the hypercomplex representation can be defined which
recovers the \emph{complex} modulating signal.
This is the \emph{complex} instantaneous amplitude.
Further, the phase of the orthonormal complex exponential may be recovered.
This is the instantaneous phase,
and its derivative is the instantaneous frequency.
In the case of modulation of a complex exponential the instantaneous
frequency is that of the complex exponential `carrier'.
Just as the classical complex analytic signal contains
both the original real signal (in the real part) and a real orthogonal signal
(in the imaginary part),
our hypercomplex signal representation contains two complex signals:
the original signal and an orthogonal signal.

We have previously published partial results on this topic
\cite{SangwineLeBihan:2007,LeBihanSangwine:2008,LeBihanSangwine:2008a}
as we developed the ideas presented here. The idea of orthonormal
complex modulation has not previously been presented.

The sequence of topics in the rest of the paper is as follows.
Section \ref{QFTsection} presents the quaternion algebra
and quaternion Fourier transform concepts
used to construct the hypercomplex representation of a complex signal
using a one-sided quaternion spectrum.
Section \ref{sec:hypercomplex} discusses the hypercomplex representation
of a complex signal. Section \ref{sec:Hmodphase} discusses the instantaneous amplitude and frequency concept for a complex signal through its hypercomplex representation. Section \ref{Examples} presents some signal examples to illustrate the ideas presented in the paper.

\section{Quaternion Fourier\\ Transform}
\label{QFTsection}

In this section, we present  the definition and
properties of the 
Quaternion Fourier transform (QFT) of complex valued
signals. Before introducing the main definitions, we give some
prerequisites on quaternions.
Interested readers may refer to \cite{Hamilton:1853,Ward1997}
for a more detailed presentation of quaternions.

\subsection{Preliminary remarks}
Quaternions were discovered by Sir W. R. Hamilton in 1843
\cite{Hamilton:1853}. Quaternions are 4D hypercomplex numbers that form a noncommutative
division ring denoted $\H$. A quaternion $q \in \H$
has a Cartesian form: $q=q_0+q_1\i+q_2\j+q_3\k$, with $q_0,q_1,q_2,q_3
\in {\mathbb R}$ and $\i,\j,\k$ roots of $-1$ satisfying
$\i^2=\j^2=\k^2=\i\j\k=-1$.
The \emph{scalar part} of $q$ is $a$: ${\cal S}(q)=q_0$.
The \emph{vector part} of $q$ is: ${\cal V}(q)=q-{\cal S}(q)$.
We will also make use of the real and imaginary parts using the following notation:
\begin{equation}
q_0=\Re(q),\quad q_1=\Im_{\i}(q),\quad q_2=\Im_{\j}(q),\quad q_3=\Im_{\k}(q)
\end{equation}
Quaternion multiplication is not commutative, so that in
general $qp\neq pq$ for $p,q \in \H$.
The conjugate of $q$ is $\qconjugate{q}=q_0-q_1\i-q_2\j-q_3\k$.
The norm of $q$ is $\norm{q}=|q|^2=(q_0^2+q_1^2+q_2^2+q_3^2)=q\qconjugate{q}$.
A quaternion with $q_0=0$ is called pure.
$|q|$ is the \emph{modulus}
of $q$. If $|q|=1$, then $q$ is called a \emph{unit} quaternion.
The inverse of $q$ is $q^{-1}=\qconjugate{q}/\norm{q}$.
Pure unit quaternions are special quaternions,
among which are $\i$, $\j$ and $\k$. Together with the identity of $\H$, they form a \emph{quaternion basis}: $\left\{1,\i,\j,\k\right\}$.

In $\H$, in addition to the conjugation,
there exist anti-involutions which play an important role, as noticed by B\"ulow \emph{et al.} \cite{BulowSommer:2001}. In the special case where $\left\{1,\i,\j,\k\right\}$ is chosen as the basis for $\H$ these involutions are, for a given quaternion $q$, the following
 \begin{equation}
 \begin{array}{ccc}
\qconjugate{q}^{\i}=-\i q \i,\,& \qconjugate{q}^{\j}=-\j q \j,\, & \qconjugate{q}^{\k}=-\k q \k.
 \end{array}
 \end{equation}
Quaternions can also be viewed as complex numbers with complex
components, \emph{i.e.} one can write $q=z_1+ z_2 \j $ in the
basis $\left\{ 1,\i,\j,\k \right\}$ with $z_1,z_2 \in \C^{\i}$, \emph{i.e.} $z_{\alpha}=\Re(z_{\alpha})+\i \Im(z_{\alpha})$ for $\alpha=1,2$. This is called the Cayley-Dickson form of $q$. 
  
Among the possible representations of $q$, we will make use in this paper of one that was recently introduced: the \emph{polar} Cayley-Dickson form \cite{10.1007/s00006-008-0128-1}. 
Any quaternion $q$ also has a unique polar Cayley-Dickson form given by: 
\begin{equation}
q = A_q \exp(B_q\j) = (a_0 + a_1 \i)\exp ((b_0+b_1\i)\j)
\label{eq:polarCD}
\end{equation}
where $A_q=a_0 + a_1 \i$ is the \emph{complex modulus} of $q$ and
$B_q=b_0+b_1\i$ its \emph{complex
phase}. As explained in \cite{10.1007/s00006-008-0128-1}, given a
\emph{degenerate} quaternion $p=(c+d\i)\j$, then its 
exponential $e^{p}$ is given by:
\begin{equation}
e^p=\cos |p| + \frac{p}{|p|} \sin |p| = \alpha + \beta \j + \delta \k
\end{equation}
leading to $\alpha = \cos |p|$, $\beta = (c/|p|) \sin |p|$ and $\delta
= (d/|p|) \sin |p|$, with $|p|=\sqrt{c^2+d^2}$. Recall also that, for
any quaternion $p$, $|e^p|=1$. Then, for any quaternion $q$ given in
its polar form as in equation (\ref{eq:polarCD}):
\begin{equation}
e^{B_q \j} = \cos|B_q|+\j\frac{b_0}{|B_q|}\sin|B_q| + \k\frac{b_1}{|B_q|}\sin|B_q|
\label{expBqj}
\end{equation}
where $|B_q| = \sqrt{b_0^2+b_1^2}$.
As a consequence, using the right hand side expression in
(\ref{eq:polarCD}), given a quaternion $q=q_0+q_1 \i + q_2 \j + q_3
\k$,  then the complex components of its polar Cayley-Dickson form,
\emph{i.e.} $A_q$ and $B_q$, are given as follows:
\begin{equation}
\left\{
\begin{array}{rcl}
A_q & = & \frac{q_0 + \i q_1}{\cos(\sqrt{q_2^2+q_3^2})}\\
B_q & = & K \left(\frac{q_0q_3+q_1q_2}{q_0^2+q_1^2} + \i \frac{q_0q_3 - q_1q_2}{q_0^2+q_1^2} \right)
\end{array}
\right.
\end{equation}
with $K=\sqrt{q_2^2+q_3^2}\arctan\left(\sqrt{q_2^2+q_3^2}\right)$. 
In Section \ref{Hmodphase}, we will
make use of the complex modulus $A_q$ for the interpretation of the
hypercomplex representation.

\subsection{Quaternion Fourier transform}
In this paper, we use a 1D version of the (right) QFT first
defined in discrete-time form in \cite{SangwineEll:1998}. 
Thus, it is necessary to specify the axis (a pure unit
quaternion) of the transform. So, we will denote by $QFT_{\mu}{\,}$ a QFT with
transformation axis $\mu$. In order to work with the classical quaternion basis $\left\{1,\i,\j,\k\right\}$ and for reasons that will be explained later, we will use $\j$ as the transform axis. The only restriction on the transform axis is that it must be orthogonal to the original basis of the signal (here $\left\{1,\i\right\}$). We now present the definition and some properties of the transform used here.

\begin{definition}
\label{def:qft}
Given a complex valued signal $z(t)=z_r(t)+\i z_i(t)$ that belongs to $L^1(\R,\C)\cap L^2(\R,\C)$, then the quaternion valued function denoted $Z_{\j}(\nu)$, and given as:
\begin{equation}
Z_{\j}(\nu)=QFT_{\j}\left[z(t)\right]=\int\limits_{-\infty}^{+\infty}z(t)e^{-\j 2 \pi \nu t} \dt
\label{Eqn:QFT2}
\end{equation}
and with $Z_{\j}(\nu) \in L^1(\R,\H)\cap L^2(\R,\H)$ is called the \emph{Quaternion Fourier Transform} $QFT_{\j}$ with respect to axis $\j$ of $z(t)$. 
\end{definition}
The inverse transform is defined as follows.
\begin{definition}
Given a quaternion Fourier transform $Z_{\j}(\nu)$ (with respect to axis $\j$), then:
\begin{equation}
z(t)=IQFT_{\j}\left[Z_{\j}(\nu)\right]=\int\limits_{-\infty}^{+\infty}Z_{\j}(\nu)e^{\j 2 \pi \nu t} d\nu
\label{Eqn:QFT3}
\end{equation}
is the \emph{Inverse Quaternion Fourier Transform} $IQFT_{\j}$ with respect to ${\j}$ of $Z_{\j}(\nu)$.
\end{definition}

\subsection{Properties of the QFT}
In order to provide insight into the advantages of using the QFT to analyze complex signals, we present a non-exhaustive list of the QFT properties. 
\begin{property} \textit{(Transform symmetry)}
\label{prop:symQFT}
Given a complex signal $z(t)=z_r(t)+\i z_i(t)$ and its quaternion
Fourier transform denoted by $Z_{\j}(\nu)$, then the following properties hold:
\begin{itemize}
\newlength{\spacer}\settowidth{\spacer}{even}
\newlength{\odd}\settowidth{\odd}{odd}
\addtolength{\spacer}{-\odd}
\item The even                part of $z_r(t)$ is in      $\Re(Z_{\j}(\nu))$
\item The odd\hspace{\spacer} part of $z_r(t)$ is in $\Im_{\j}(Z_{\j}(\nu))$
\item The even                part of $z_i(t)$ is in $\Im_{\i}(Z_{\j}(\nu))$
\item The odd\hspace{\spacer} part of $z_i(t)$ is in $\Im_{\k}(Z_{\j}(\nu))$
\end{itemize}
\end{property}
\begin{proof}
Expand \eqref{Eqn:QFT2} into real and imaginary parts with respect to $\i$,
and expand the quaternion exponential into cosine and sine components:
\begin{align*}
Z_{\j}(\nu)
&=  \!\int\limits_{-\infty}^{+\infty}{\!\left[z_r(t)+\i z_i(t)\right]}
            \left[\cos(2\pi\nu t) - \j\sin(2\pi\nu t)\right]\dt\\ 
&=  \!\int\limits_{-\infty}^{+\infty}\!z_r(t)\cos(2\pi\nu t)\dt 
 -\j\!\int\limits_{-\infty}^{+\infty}\!z_r(t)\sin(2\pi\nu t)\dt\\
&+\i\!\int\limits_{-\infty}^{+\infty}\!z_i(t)\cos(2\pi\nu t)\dt 
 -\k\!\int\limits_{-\infty}^{+\infty}\!z_i(t)\sin(2\pi\nu t)\dt 
\end{align*}
from which the stated properties are evident.
\end{proof}
It is a known fact that the \emph{classical} Fourier transform of a real signal possesses Hermitian symmetry. It is also known that the \emph{classical} Fourier transform of a complex signal possesses no symmetry at all. The QFT of a complex signal has symmetry according
to the following properties: 
\begin{property}  \textit{($\i$-involution reversal)}
\label{prop:hermq}
Given a complex valued signal $z(t)$ and its quaternion
Fourier transform denoted $QFT_{\j}[z(t)]=Z_{\j}(\nu)$, then the following property holds:
\begin{equation}
Z_{\j}(-\nu)=-\i Z_{\j}(\nu) \i = \qconjugate{Z_{\j}(\nu)}^{\i}
\end{equation}
\end{property}
\begin{proof}
This can be directly checked from Property~\ref{prop:symQFT}.
\end{proof}

% Real and i-imaginary parts
\begin{figure}
\setlength{\unitlength}{0.01\columnwidth}
\begin{center}
\subfigure{
\begin{picture}(45,50)
%\put(0,0){\dashbox(45,50){}}
\thicklines
\put(0,20){\vector(1,0){40}} % X-axis
\put(20,0){\vector(0,1){40}} % Y-axis
\put(15,42){\makebox(10,5){$\Re(Z_{\j}(\nu))$}} % Y-axis Label
\put(40,18){\makebox(6,4){$\nu$}} % Y-axis Label
\put(16,15){\makebox(4,4){$0$}} % O position
\put(8,20){\line(1,2){4}}
\put(12,28){\line(1,0){8}}
\put(20,28){\line(1,0){8}}
\put(28,28){\line(1,-2){4}}
\end{picture}
}
\subfigure{
\begin{picture}(45,50)
%\put(0,0){\dashbox(45,50){}}
\thicklines
\put(0,20){\vector(1,0){40}} % X-axis
\put(20,0){\vector(0,1){40}} % Y-axis
\put(15,42){\makebox(10,5){$\Im_{\i}(Z_{\j}(\nu))$}} % Y-axis Label
\put(40,18){\makebox(6,4){$\nu$}} % Y-axis Label
\put(16,15){\makebox(4,4){$0$}} % O position
\put(8,20){\line(1,2){4}}
\put(12,28){\line(1,0){8}}
\put(20,28){\line(1,0){8}}
\put(28,28){\line(1,-2){4}}
\end{picture}
}
% j-imaginary and k-imaginary parts
\subfigure{
\begin{picture}(45,50)
%\put(0,0){\dashbox(45,50){}}
\thicklines
\put(0,20){\vector(1,0){40}} % X-axis
\put(20,0){\vector(0,1){40}} % Y-axis
\put(15,42){\makebox(10,5){$\Im_{\j}(Z_{\j}(\nu))$}} % Y-axis Label
\put(40,18){\makebox(6,4){$\nu$}} % Y-axis Label
\put(16,15){\makebox(4,4){$0$}} % O position
\put(8,20){\line(1,-2){4}}
\put(12,12){\line(1,0){8}}
\put(20,28){\line(1,0){8}}
\put(28,28){\line(1,-2){4}}
\end{picture}
}
\subfigure{
\begin{picture}(45,50)
%\put(0,0){\dashbox(45,50){}}
\thicklines
\put(0,20){\vector(1,0){40}} % X-axis
\put(20,0){\vector(0,1){40}} % Y-axis
\put(15,42){\makebox(10,5){$\Im_{\k}(Z_{\j}(\nu))$}} % Y-axis Label
\put(40,18){\makebox(6,4){$\nu$}} % Y-axis Label
\put(16,15){\makebox(4,4){$0$}} % O position
\put(8,20){\line(1,-2){4}}
\put(12,12){\line(1,0){8}}
\put(20,28){\line(1,0){8}}
\put(28,28){\line(1,-2){4}}
\end{picture}
}
\end{center}
\caption{\label{draw:symm}Symmetry of the four components of the quaternion valued Fourier transform $Z_{\j}(\nu)$ of a complex signal $z(t)$.}
\end{figure}
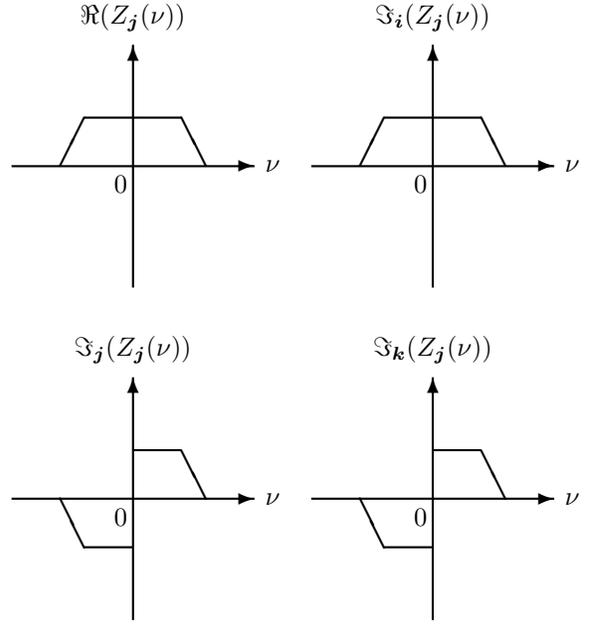
This symmetry property of the QFT of a complex signal is illustrated in
\figurename~\ref{draw:symm}
and is central to the justification of using the QFT to
analyze a complex-valued signal carrying complementary but different
information in its real and imaginary parts. Using the QFT, it is
possible to keep the odd and even parts of the real and imaginary parts of
the signal in four different components in the transform
domain. This idea was also the initial motivation of B\"ulow,
Sommer and Felsberg when they developed the monogenic signal for images
\cite{Bulow:1999,BulowSommer:2001,FelsbergSommer:2001}. Note that the
use of hypercomplex Fourier transforms was
originally introduced in 2D Nuclear Magnetic Resonance image analysis \cite{Ernst:1987,Delsuc:1988}.

For use later in the paper, we mention another property.
\begin{property} \textit{($\j$-involution conjugation)}
\label{prop:qftconj}
Given a complex signal $z(t)$ and its $QFT_{\j}$ denoted $Z_{\j}(\nu)$, then the $QFT_{\j}$ of its conjugate $\cconjugate{z}(t)$ is given by:
\begin{equation}
QFT_{\j}\left[\cconjugate{z}(t)\right]=-\j Z_{\j}(\nu) \j = \qconjugate{Z_{\j}(\nu)}^{\j}
\end{equation}
\end{property}
\begin{proof}
This can be directly checked by calculation.	
\end{proof}
An important issue for the forthcoming study of complex signals is the way \emph{modulation} affects the QFT.
Of special importance is \emph{ortho-complex} modulation.
\begin{property} \textit{(Ortho-complex modulation)}
\label{prop:orthoC}
Given a complex signal $z(t)$ and its $QFT_{\j}$ denoted $Z_{\j}(\nu)$, the $QFT_{\j}$ of the ortho-complex modulation of $z(t)$ by the quaternion-valued exponential $\exp(\j 2\pi \nu_0 t)$ is
\begin{equation}
QFT_{\j}[z(t)e^{\j 2\pi \nu_0 t}] = Z_{\j}(\nu-\nu_0)
\label{eq:freqshift}
\end{equation}
\end{property}
\begin{proof}
As $e^{\j 2\pi \nu_0 t}$ is an exponential with the same axis $\j$ as the $QFT_{\j}$, then direct calculation similar to the complex case can be performed as $e^{\j a}e^{\j b}=e^{\j(a+b)}$. The desired result is then directly obtained.
\end{proof}
Note that a similar frequency shift occurs simultaneously in the four components of the quaternion valued Fourier transform.
Note also that the axis of the modulation plays an important role as a modulation using axis $\i$ or $\k$ would not lead to a simple shift of the spectrum\footnote{This is due to the fact that for two real numbers $a$ and $b$, then $e^{\j a}e^{\k b} \neq e^{\j a +\k b}$. This is in turn due to the fact that $\j$ and $\k$ do not commute. This is a special case of the well known
Baker-Campbell-Hausdorff formula \cite{Chirikjian2000}.}. Also, it is obvious from this property and the convolution property (see equation \ref{Eqn:QFT3A})
fulfilled by the $QFT_{\j}$ that $QFT_{\j}[e^{\j 2\pi \nu_0 t}] = \delta(\nu-\nu_0)$.

% Real and i-imaginary parts
\begin{figure}[ht!]
\setlength{\unitlength}{0.01\columnwidth}
\begin{center}
\subfigure{
\begin{picture}(45,50)
%\put(0,0){\dashbox(45,50){}}
\thicklines
\put(0,20){\vector(1,0){40}} % X-axis
\put(20,0){\vector(0,1){40}} % Y-axis
\put(15,42){\makebox(10,5){$\Re(Z_{\j}(\nu-\nu_0))$}} % Y-axis Label
\put(40,18){\makebox(6,4){$\nu$}} % Y-axis Label
\put(16,15){\makebox(4,4){$0$}} % O position
\put(27.5,16){\makebox(4,4){$\nu_0$}} % nu_O position
\put(27,19.5){\line(0,1){1}}
\put(22,20){\line(1,2){4}}
\put(26,28){\line(1,0){1}}
\put(27,28){\line(1,0){1}}
\put(28,28){\line(1,-2){4}}
\end{picture}
}
\subfigure{
\begin{picture}(45,50)
%\put(0,0){\dashbox(45,50){}}
\thicklines
\put(0,20){\vector(1,0){40}} % X-axis
\put(20,0){\vector(0,1){40}} % Y-axis
\put(15,42){\makebox(10,5){$\Im_{\i}(Z_{\j}(\nu-\nu_0))$}} % Y-axis Label
\put(40,18){\makebox(6,4){$\nu$}} % Y-axis Label
\put(16,15){\makebox(4,4){$0$}} % O position
\put(27.5,16){\makebox(4,4){$\nu_0$}} % nu_O position
\put(27,19.5){\line(0,1){1}}
\put(22,20){\line(1,2){4}}
\put(26,28){\line(1,0){1}}
\put(27,28){\line(1,0){1}}
\put(28,28){\line(1,-2){4}}
\end{picture}
}
% j-imaginary and k-imaginary parts
\subfigure{
\begin{picture}(45,50)
%\put(0,0){\dashbox(45,50){}}
\thicklines
\put(0,20){\vector(1,0){40}} % X-axis
\put(20,0){\vector(0,1){40}} % Y-axis
\put(15,42){\makebox(10,5){$\Im_{\j}(Z_{\j}(\nu-\nu_0))$}} % Y-axis Label
\put(40,18){\makebox(6,4){$\nu$}} % Y-axis Label
\put(16,15){\makebox(4,4){$0$}} % O position
\put(27.5,16){\makebox(4,4){$\nu_0$}} % nu_O position
\put(22,20){\line(1,-2){4}}
\put(26,12){\line(1,0){1}}
\put(27,12){\line(0,1){16}}
\put(27,28){\line(1,0){1}}
\put(28,28){\line(1,-2){4}}
\end{picture}
}
\subfigure{
\begin{picture}(45,50)
%\put(0,0){\dashbox(45,50){}}
\thicklines
\put(0,20){\vector(1,0){40}} % X-axis
\put(20,0){\vector(0,1){40}} % Y-axis
\put(15,42){\makebox(10,5){$\Im_{\k}(Z_{\j}(\nu-\nu_0))$}} % Y-axis Label
\put(40,18){\makebox(6,4){$\nu$}} % Y-axis Label
\put(16,15){\makebox(4,4){$0$}} % O position
\put(27.5,16){\makebox(4,4){$\nu_0$}} % nu_O position
\put(22,20){\line(1,-2){4}}
\put(26,12){\line(1,0){1}}
\put(27,12){\line(0,1){16}}
\put(27,28){\line(1,0){1}}
\put(28,28){\line(1,-2){4}}
\end{picture}
}
\end{center}
\caption{\label{draw:symm2}Four components of the quaternion valued Fourier transform $Z_{\j}(\nu-\nu_0)$ of the orthocomplex modulated signal $z(t)e^{\j 2 \pi \nu_0 t}$ and where the complex signal $z(t)$ has a $QFT_{\j}$ verifying $\left|Z_{\j}(\nu)\right|=0$ for $\left|\nu\right|>\nu_0$. This spectrum is an example of \eqref{eq:freqshift}.}  
\end{figure}
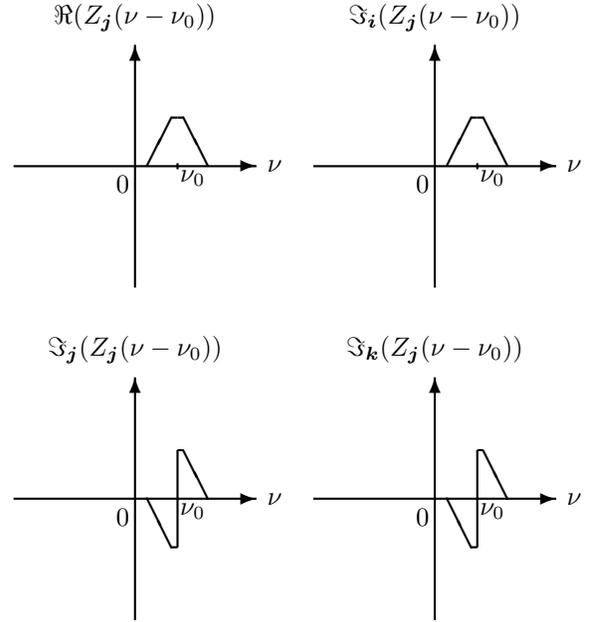

\begin{property} \textit{(Positive spectra modulation)}
Given a complex signal $z(t)$ with a band-limited $QFT_{\j}$: $X_{\j}(\nu)=0$ for $|\nu|>\nu_m$ for a given $\nu_m$. Then, for some $\nu_0 > \nu_m$, the ortho-complex modulated signal $z(t)e^{\j 2 \pi \nu_0 t}$ has a right-sided $QFT_{\j}$.
\end{property}
\begin{proof}
This results follows directly from property \ref{prop:orthoC} and the band-limited property of $z(t)$.
\end{proof}
This property is illustrated in \figurename~\ref{draw:symm2}.

Now, we introduce the convolution property in the case where one of
the signals is real valued. Note that we do not give the general case as
only the above mentioned one will be of use in the sequel. 

\begin{property} \textit{(Convolution)}
Given $g$ and $f$ such that: $g:{\mathbb R}\rightarrow{\mathbb C}$ and $f:{\mathbb R}\rightarrow{\mathbb
R}$, the $\QFT_{\j}$ of their convolution product $(g*f)(t)$ is:
\begin{equation}
QFT_{\j}\left[(g*f)(t)\right]   = QFT_{\j}\left[g(t)\right] \, QFT_{\j}\left[f(t)\right]
\label{Eqn:QFT3A}
\end{equation}
\end{property}
\begin{proof}
Taking the $\QFT_{\j}$ of their convolution product one gets the result
shown in \eqref{eqn:wide}
\begin{figure*} % This is a two-column figure without caption to contain the equation
                % which will not fit into one column.
\begin{equation}
\label{eqn:wide}
\begin{aligned}
QFT_{\j}\left[(g*f)(t)\right] &= \int\limits_{-\infty}^{+\infty}
                     \int\limits_{-\infty}^{+\infty}
                     g(\tau)f(t-\tau)\dtau\,
                     e^{-2\j \pi \nu t} \dt\\
                  % &= \int\limits_{-\infty}^{+\infty}
                  %    \int\limits_{-\infty}^{+\infty}
                  %    g(\tau)e^{-\j 2 \pi \nu (t'+\tau)}f(t')d\tau \dt'\\
                  &= \int\limits_{-\infty}^{+\infty}
                     g(\tau)e^{-2\j \pi \nu \tau}\dtau
                     \int\limits_{-\infty}^{+\infty}
                    f(t')e^{-2\j \pi \nu t'}\dt'\\
                  &= QFT_{\j}\left[g(t)\right]QFT_{\j}\left[f(t)\right]
              %    &= QFT_{\j}\left[g(t)\right]\qconjugate{QFT_{\j}\left[f(t)\right]}^{\j}
\end{aligned}
\end{equation}
\rule{0.99\textwidth}{0.5pt}
\end{figure*}
which is true thanks to the fact that $f(t)$ commutes with the $\j$
exponential as $f(t)$ is real valued.
\end{proof}
Note that the order matters in the product of the $QFT_{\j}$s in
\eqref{Eqn:QFT3A} as they are quaternion and complex valued (which do
not commute in the general case). This expression is thus valid only
when $f$ is convolved on the right. Convolution on the left would result
in the appearance of a conjugation in the Fourier domain.

\subsection{Hilbert transform}
\label{subsec:ConvHilb}
The quaternion Fourier transform
of the \emph{Hilbert transformer} $x(t)=1/\pi t$ is $-\j\sign(\nu)$,
where $\j$ is the axis of the transform. This arises simply from the fact that using the exponential with axis $\j$ in the transform will lead to a real and $\Im_{\j}$ part for the $QFT_{\j}$ only, which is a degenerate quaternion isomorphic to a complex number. Substituting $x(t)$ into \eqref{Eqn:QFT2}, we get:
\begin{equation*}
QFT_{\j}\left[\frac{1}{\pi t}\right]=\frac{1}{\pi}\int\limits_{-\infty}^{+\infty}\frac{e^{-\j 2\pi\nu t}}{t}\,\dt
\end{equation*}
and this is clearly isomorphic to the classical complex case. The solution
in the classical case is $-{\i}\sign(\nu)$, and hence in the quaternion case
must be as stated above.

It is straightforward to see that, given an arbitrary real signal $y(t)$,
subject only to the constraint that its classical Hilbert transform
$\hilbert{y(t)}$ exists, then one can easily show that the classical Hilbert transform
of the signal may be obtained using a quaternion Fourier transform as
follows:
\begin{equation}
\label{eqn:hilbertbyqft}
\hilbert{y(t)} = IQFT_{\j}\left[-\j\sign(\nu)Y_{\j}(\nu)\right]
\end{equation}
where $Y_{\j}(\nu)=QFT_{\j}\left[y(t)\right]$. This result follows from the isomorphism between the quaternion and complex Fourier transforms when operating on a real signal, and it may be seen to be
the result of a convolution between the signal $y(t)$ and the quaternion
Fourier transform of $x(t)=1/\pi t$.
Note that $\j$ and $Y_{\j}(\nu)$ commute as a consequence of $y(t)$ being
real. We can now define the Hilbert transform of a \emph{complex} signal.

\begin{definition}
\label{def:hyperhilbert}
Consider a complex signal $z(t)=z_r(t)+\i z_i(t)$
and its quaternion Fourier transform $Z_{\j}(\nu)$ as defined in Definition~\ref{def:qft}.
Then, the hypercomplex analogue of the Hilbert transform of $z(t)$, is as follows:
\begin{equation}
\qhilbert{\j}{z(t)} =\iqfourier{\j}{-\j\sign(\nu)Z_{\j}}(\nu)
\end{equation}
where the Hilbert transform is thought of as:
$\qhilbert{\j}{z(t)}=p.v.\left(z*\frac{1}{\pi t}\right)$, where the
principal value (p.v.) is understood in its classical way. This result follows from \eqref{eqn:hilbertbyqft} and the linearity
of the quaternion Fourier transform.
To extract the imaginary part, the vector part of the quaternion signal must
be multiplied by $-\i$.
An alternative is to take the scalar or inner product of the vector part with
$\i$.
Note that $\j$ and $Z_{\j}(\nu)$ anticommute because $\j$ is orthogonal
to $\i$, the axis of $Z_{\j}(\nu)$.
Therefore the ordering is not arbitrary, but changing it simply conjugates
the result.
\end{definition}
\emph{Note that the signal $z(t)$ is considered to be non-analytic in the classical (complex) sense,
that is its real and imaginary parts are \emph{not} orthogonal. 
However, this definition is valid if $z(t)$ is analytic, as it can be considered as a degenerate case of the more
general non-analytic case.} We end this section with a useful property of the right-sided $QFT_{\j}$. 

\begin{property} \textit{(Positive spectrum quaternion signal)}
Consider a $QFT_{\j}$ denoted by $Z_{\j}(\nu)$. Then, $Z_{\j}(\nu)$ is said to be right-sided if $Z_{\j}(\nu)=0$, $\forall \nu < 0$. Now, if a $QFT_{\j}$ denoted $Z_{\j}(\nu)$ is right-sided, then its Inverse Quaternion Fourier transform $IQFT_{\j}\left[Z_{\j}(\nu)\right]$ is a quaternion valued signal.
\end{property}
\begin{proof}
The $QFT_{\j}$ of a complex signal has the symmetry described in Property \ref{prop:hermq} and illustrated in \figurename~\ref{draw:symm}. Assuming these symmetries, a right-sided $QFT_{\j}$ can only be obtained from a $QFT_{\j}$, denoted $Z_{\j}(\nu)$, by the following linear combination: $(1+\sign(\nu))Z_{\j}(\nu)$ which simply cancels out the negative frequencies. Now, the $IQFT_{\j}$ of $\sign(\nu)$ is, as previously seen, equal to $\j/\pi t$. The convolution of $1/\pi t$ with $z(t)$ leads to a complex signal with
$\Re$ and $\Im_{\i}$
which, when multiplied by $\j$ becomes a degenerate quaternion signal with only the $\Im_{\j}$ and $\Im_{\k}$ parts non-zero. Now, the total $IQFT_{\j}$ of $(1+\sign(\nu))Z_{\j}(\nu)$ consists of $z(t)+{\j}\left(\frac{1}{\pi t}*z(t)\right)(t)$ which is thus made of a real part and three imaginary parts of axis ${\i}$, ${\j}$ and ${\k}$. 
\end{proof}

In the next section, we make use of the QFT and Hilbert transform to build up a hypercomplex representation for complex valued signals.

\section{The hypercomplex signal\\ representation}
\label{sec:hypercomplex}

We define a hypercomplex representation, denoted $\hat{z}(t)$, of the complex signal $z(t)$
by a similar approach to that originally developed by Ville \cite{Ville:1948}. The $\hat{z}(t)$ signal is quaternion-valued. The following definitions give the details of the construction of this signal.
\begin{property} \textit{(Hypercomplex extension)}
Given a complex signal $z(t)$, one can associate to it a unique
canonical pair corresponding to a (complex) modulus and phase.
The modulus and phase are uniquely defined through the hypercomplex
representation of the complex signal, which is quaternion valued.
\end{property}
\begin{proof}
Cancelling the negative frequencies of the QFT leads to a quaternion
signal in the time domain. Then, any quaternion signal has a modulus
and phase defined using its CD polar form.
\end{proof}
\begin{definition}
\label{def:hyperanalytic}
Given a complex valued signal $z(t)$ that can be expressed as $z(t)=z_r(t)+\i z_i(t)$,
then the hypercomplex representation of $z(t)$,
denoted $\hat{z}(t)$ is given by: 
\begin{equation}
\hat{z}(t)=z(t)+\j\qhilbert{\j}{z(t)}
\label{defHana}
\end{equation}
where $\qhilbert{\j}{z(t)}$ is the hypercomplex analogue of the
Hilbert transform of $z(t)$ defined in the preceding definition.
The $QFT_{\j}$ of this hypercomplex representation of the complex signal,
denoted $\hat{Z}_{\j}(\nu)$, is thus:
\begin{align*}
\hat{Z}_{\j}(\nu)
&=Z_{\j}(\nu)-\j^{2}\sign(\nu)Z_{\j}(\nu)\\
&=\left[1+\sign(\nu)\right]Z_{\j}(\nu)\\
& = 2 U(\nu) Z_{\j}(\nu)
\end{align*}
where $U(\nu)$ is the classical unit step function.
\end{definition}
This result is unique to the quaternion Fourier transform representation of the
hypercomplex representation, which has a one-sided quaternion Fourier spectrum.
This means that the hypercomplex representation may be constructed from a complex signal
$z(t)$ in exactly the same way that an analytic signal may be constructed from a
real signal, by suppression of negative frequencies in the Fourier domain.
The only difference is that in the hypercomplex representation,
a quaternion rather than a complex Fourier transform must be used,
and of course the complex signal must be put in the form $z(t)=z_r(t)+\i z_i(t)$
which, although a quaternion signal, is isomorphic to the original complex signal.

\begin{property}  \textit{(Hypercomplex phasor form)}
Given a complex signal $z(t)$, one can associate to it a unique
canonical pair $\left[\rho(t),\phi(t)\right]$ consisting of a modulus and a phase.
The modulus and phase are uniquely defined through the hypercomplex
representation of the complex signal.
\end{property}
\begin{proof}
Cancelling the negative frequencies of the $QFT_{\j}$ of a complex signal leads to a quaternion
signal in the time domain, with similar frequency content (the original $QFT_{\j}$ of the complex signal can be recovered thanks to the symmetry properties of the quaternion Fourier transform).
As any quaternion signal has a unique modulus
and phase defined through its CD polar form, those are uniquely associated to the original complex signal.
\end{proof}

A second important property of the hypercomplex representation of a complex signal
$z(t)$ is that it maintains a separation between the different even and odd parts of
the original signal.

\begin{property}  \textit{(Orthogonal split)}
The original signal $z(t)$ is the \emph{simplex} part \cite[Theorem 1]{10.1109/ICIP.2000.899828},
\cite[\S\,13.1.3]{Multivariate:2009} of its corresponding hypercomplex representation $\hat{z}(t)$.
The \emph{perplex} part is the orthogonal or `quadrature' component, $o(t)$. 
They are obtained by:
\begin{align}
\label{eqn:simplex}
z(t)&=\frac{1}{2}\left(\hat{z}(t) + \qconjugate{\hat{z}(t)}^{\i}\right)\\
\label{eqn:perplex}
o(t)&=\frac{1}{2}\left(z_{+}(t) - \qconjugate{\hat{z}(t)}^{\i}\right)
\end{align}
\end{property}
\begin{proof}
This follows from \eqref{defHana}.
Writing this in full by substituting the orthogonal signal for $\qhilbert{\j}{z(t)}$:
\[
\hat{z}(t) = z(t) + \j o(t) = z_r(t) + \i z_i(t) + \j o_r(t) -\k o_i(t) 
\]
and substituting this into equation \eqref{eqn:simplex}, we get:
\begin{align*}
z(t)&=\frac{1}{2}\left(
\begin{aligned}&\phantom{\i\left[\right.}z_r(t) + \i z_i(t) + \j o_r(t) - \k o_i(t)\\
              -&         \i\left[        z_r(t) + \i z_i(t) + \j o_r(t) - \k o_i(t) \right]\i
\end{aligned}\right)\\
\intertext{and since $\i$ and $\j$ are orthogonal unit pure quaternions, $\i\j=-\j\i$:} 
    &=\frac{1}{2}\left(
\begin{aligned}&\phantom{\i\left[\right.}z_r(t) + \i z_i(t) + \j o_r(t) - \k o_i(t)\\
              +&\phantom{\i\left[\right.}z_r(t) + \i z_i(t) - \j o_r(t) + \k o_i(t)
\end{aligned}\right)
\end{align*}
from which the first part of the result follows. Equation \eqref{eqn:perplex} differs
only in the sign of the second term, and it is straightforward to see that if $\hat{z}(t)$
is substituted, $z(t)$ cancels out, leaving $o(t)$. 
\end{proof}

\section{Instantaneous amplitude,\\ phase and frequency}
\label{sec:Hmodphase}

\label{Hmodphase}
Using the hypercomplex representation $\hat{z}(t)$ of a complex signal $z(t)$,
we now investigate the concept of instantaneous amplitude and phase for $z(t)$.

Recall that $\hat{z}(t)$ is quaternion valued,
and thus it can be expressed using its \emph{polar} Cayley-Dickson form
given in
\eqref{eq:polarCD} as:
\begin{equation}
\hat{z}(t)=\hat{\rho}(t)e^{\hat{\phi}(t) \j}
\end{equation}
where $\hat{\rho}(t)$ is called the \emph{instantaneous amplitude} of
$\hat{z}(t)$ and $\hat{\phi}(t)$ is the \emph{instantaneous phase} of
$\hat{z}(t)$. The pair $\left[\hat{\rho}(t),\hat{\phi}(t)\right]$ is
\emph{canonical} and associated to $z(t)$. The one-to-one
correspondence between $z(t)$ and the canonical pair can be
demonstrated using the same arguments used by Picinbono in
\cite{Picinbono:1997}. Indeed, as presented in the previous sections,
$\hat{z}(t)$ has a one-sided spectrum. This ensures the spectral
characterization (see \cite[II.B]{Picinbono:1997}) of the canonical pair.
In addition,
the one-to-one correspondence is ensured because $e^{\j
  \hat{\phi}(t)}$ is a \emph{phase signal} with respect to the
$QFT_{\j}$. Now, just as explained in \cite[III]{Picinbono:1997}, $e^{\j
  \hat{\phi}(t)}$ is a \emph{unimodular} signal, and in order for it to
be one-sided, some restrictions\footnote{Note
that the considerations given in \cite{Picinbono:1997} for
instantaneous phase, amplitude and frequency and their link with the
analytic (one-sided) signal are based on the Bedrosian theorem which
extends easily to the $QFT_{\j}$ case due to the properties
presented in Section \ref{QFTsection} and \ref{sec:hypercomplex}.} apply to $\hat{\phi}(t)$. In
short, the phase $\hat{\phi}(t)$ is of the form:
\begin{equation}
\hat{\phi}(t)=\theta+2\pi\nu_0 t + \hat{\phi}_b(t)
\end{equation}
where $\theta$ is a constant, $\nu_0$ is a given frequency and
$\hat{\phi}_b(t)$  is a band-limited contribution with maximum
frequency $B>\nu_0$. In Section \ref{Examples}, we present examples of the
former two terms.

With all these considerations, one can state the following. 
\begin{definition}
Given a complex signal $z(t)$, its instantaneous amplitude
$\hat{\rho}(t)$ and instantaneous phase $\hat{\phi}(t)$ are obtained
through the \emph{polar} Cayley-Dickson form of its hypercomplex
representation $\hat{z}(t)$. 
\end{definition}
As noted above, through this definition,
$\left[\hat{\rho}(t),\hat{\phi}(t)\right]$ is a canonical pair
associated to $z(t)$. Note that, unlike in the \emph{classical} case of
the analytic signal associated to a real signal, the \emph{instantaneous amplitude}
of $z(t)$ is complex valued.
However, it represents the low-frequency part (non modulated) of $z(t)$.
The instantaneous frequency can now be defined similar to the \emph{classical} case.
\begin{definition} 
The \emph{instantaneous frequency} of the complex signal $z(t)$ is defined as:
\begin{equation}
\hat{f}(t)=\frac{1}{2\pi}\frac{d\hat{\phi}(t)}{dt}
\end{equation}
where $\hat{\phi}(t)$ is the instantaneous phase of $z(t)$.
\end{definition} Like the classical case, this definition follows naturally from the Taylor series for $\hat{\phi}(t)$ about $t_0$, namely, 
\begin{align*}
\hat{\phi}(t)&=\hat{\phi}(t_0) + (t-t_0)\frac{d\hat{\phi}(t_0)}{dt} + R\\
                  &= ( \hat{\phi}(t_0) -t_0\frac{d\hat{\phi}(t_0)}{dt} ) + \frac{d\hat{\phi}(t_0)}{dt} t + R
\end{align*}
where $R$ is small when $t$ is close to $t_0$ and we see that $\frac{d\hat{\phi}(t_0)}{dt}$ has the role of frequency. Unlike the classical case, the complex nature of $z(t)$ requires that the oscillation direction
evolves over time (as can be seen for example in \figurename~\ref{modulation_1a} and
\ref{modulation_3a}).
Information on the orientation of this oscillation,
\emph{i.e.} the osculating plane orientation,
can be recovered from the instantaneous amplitude $\hat{\rho}(t)$.
More precisely, the normal to the osculating plane, denoted $\hat{n}(t)$,  is: 
\begin{equation}
\hat{n}(t)=\hat{\rho}(t) \times \hat{\rho}'(t)
\label{eq:hatn}
\end{equation}
where $\hat{\rho}^\prime(t)=d\hat{\rho}(t)/dt$ and where $\times$ is
the cross product of vectors corresponding to the time-varying complex
signals: $[z_r(t)~ , ~ z_i(t) ~ , ~t]$ is the vector associated
with the complex signal $z(t)=z_r(t)+\i z_i(t)$. From a geometric
point of view, $\hat{\rho}^\prime(t)$ can
be considered as the \emph{velocity} of $z(t)$.

Such a \emph{geometrical} description of a complex signal can be
achieved with the use of concepts defined through the $QFT_{\j}$. Such
an analysis could be developed to higher dimensional signals using
geometric algebra Fourier transforms.

\section{Examples}
\label{Examples}

We present three examples in which a band-limited complex signal
$A(t)$ modulates an orthogonal complex exponential $\exp(B(t)\j)$.
This can be represented using the polar Cayley-Dickson form of
\eqref{eq:polarCD} as:
\begin{equation}
q(t) = A(t) \exp(B(t)\j)
\end{equation}
The result is a quaternion-valued signal.
We can extract from this, two complex signals using the Cartesian
Cayley-Dickson form: $q(t) = z(t) + o(t)\j$, where $z(t) = z_r(t) + z_i(t)\i$
and similarly for $o(t)$. In what follows we take the first of these to be our modulated
complex signal (the second is a type of `quadrature' signal).

In the first example, the orthonormal complex exponential has constant
frequency, so that simply:
\begin{equation}
B_1(t) = 2\pi\nu_0 t
\end{equation}
 where $\nu_0$ is an arbitrary constant value. 
This example can be generated in the frequency domain,
since the modulation in the time domain corresponds to a shift in
the frequency domain by $\nu_0$. The complex signal is $z_1(t)$.

In the second example, $B_2(t)$ is a step function of the form:
\begin{equation}
B_2(t)=
\left\{
\begin{array}{lcr}
2\pi\nu_0 t  & \text{if} & \hfill t \in [0,0.25] \\
 
2\pi \nu_1 t & \text{if} & \hfill t \in [0.25,0.75] \\

2 \pi \nu_0 t & \text{if} & \hfill t \in [0.75,1]

\end{array}
\right.
\label{eq:example2_phase}
\end{equation}
where $\nu_0$ and $\nu_1$ are arbitrary constant values. This example
is generated in the same manner as the first one. The complex signal
is thus $z_2(t)$ for this example.

In the third example, the frequency of the complex exponential sweeps
linearly from an initial value to double the initial value and back to
the initial value. That is:
\begin{equation}
B_3(t)=\alpha(\Pi_{T/2}*\Pi_{T/2})(t-T/2)
\end{equation} 
where $\Pi$ is the boxcar function with $T=0.5$ and $\alpha$ a constant in the present case.
This example is also generated in the time domain. The complex signal
is denoted $z_3(t)$.

We present in the sequel how the modulating complex signal,
and the frequency of the carrier can be recovered from the quaternion
hypercomplex representation of the complex signal using the instantaneous
amplitude (the `envelope') and the instantaneous phase, both of which are
computed from the Cayley-Dickson polar form of the quaternion-valued hypercomplex
representation of the complex signals $z_1(t)$, $z_2(t)$ and $z_3(t)$.

\begin{figure}[t]
\includegraphics[width=\columnwidth]{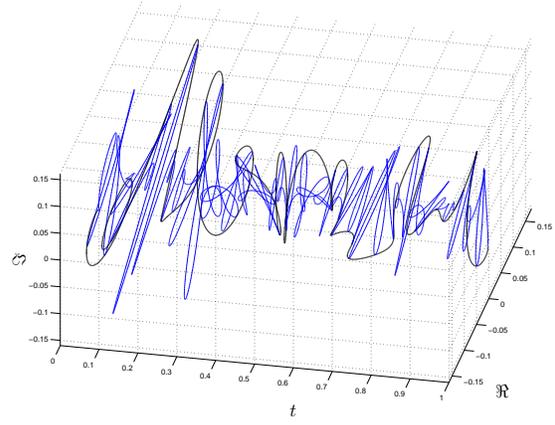}
\caption{\label{modulation_1a}Band-limited signal $z_1(t)$ (blue) and
  envelope $\hat{\rho}_1(t)$ (black) of the hypercomplex
  representation $\hat{z}_1(t)$ associated to $z_1(t)$.}
\end{figure}

%%% Example 1 
\figurename~\ref{modulation_1a} shows a modulated complex signal
$z_1(t)$ (in blue) constructed by
modulating an orthonormal complex exponential with a band-limited complex signal.
The latter was constructed from a pseudo-random signal by filtering in the
frequency domain to suppress all frequencies higher than an upper limit of 16
cycles over the length of the signal.
The orthocomplex exponential in this case has constant frequency
four times that of the highest frequency of the baseband complex signal,
and the modulated signal was created by frequency shifting the quaternion
Fourier transform of the baseband signal, as in
\eqref{eq:freqshift}. The complex envelope $\hat{\rho}_1(t)$  of the hypercomplex
representation of $z_1(t)$, {\em i.e.} $\hat{z}_1(t)$, is displayed in black in
\figurename~\ref{modulation_1a}. 

\begin{figure}
\includegraphics[width=\columnwidth,height=5cm]{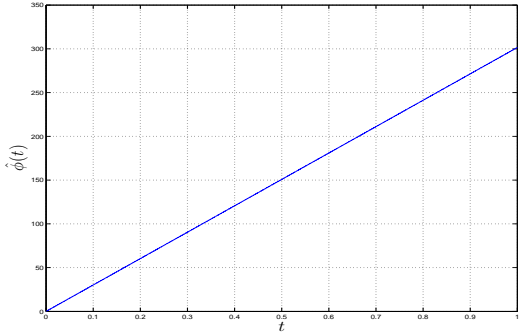}
\caption{\label{modulation_1b}Phase $\hat{\phi}_1(t)$ recovered from the polar Cayley-Dickson
form of the hypercomplex representation $\hat{z}_1(t)$ of  the signal
$z_1(t)$ displayed in \figurename~\ref{modulation_1a},
showing the linearly increasing phase $B_1(t)$.}
\end{figure}

\figurename~\ref{modulation_1b} shows the instantaneous phase %frequency
$\hat{\phi}_1(t)$ of the hypercomplex representation $\hat{z}_1(t)$ of
the signal $z_1(t)$ displayed in \figurename~\ref{modulation_1a}. This
phase is obtained using the polar Cayley-Dickson phase of the
hypercomplex representation $\hat{z}_1(t)$.

\begin{figure}[t!]
\includegraphics[width=\columnwidth]{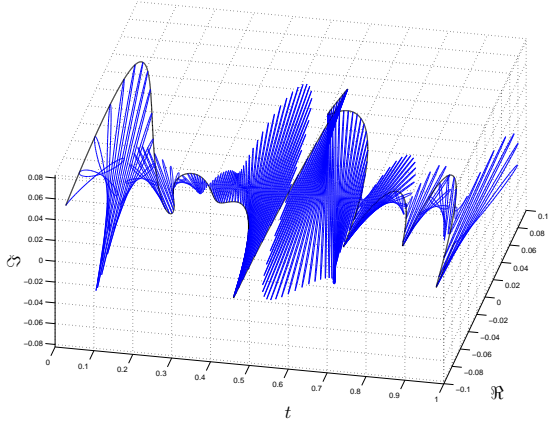}
\caption{\label{modulation_2a}Band-limited signal $z_2(t)$ (blue)  envelope $\hat{\rho}_2(t)$ (black) of the hypercomplex
  representation $\hat{z}_2(t)$ associated to $z_2(t)$.}
\end{figure}

\begin{figure}[ht]
\includegraphics[width=\columnwidth,height=5cm]{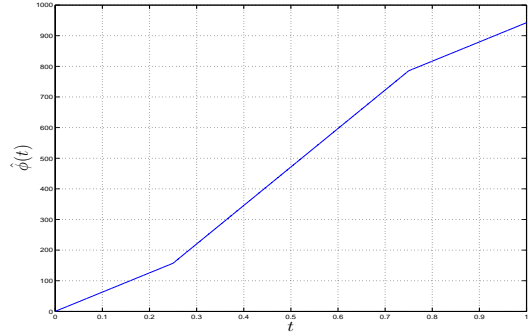}
\caption{\label{modulation_2b}Phase $\hat{\phi}_2(t)$  recovered from the polar Cayley-Dickson
                              form of the hypercomplex
                              representation $\hat{z}_2(t)$ of the
                              signal $z_2(t)$ in \figurename~\ref{modulation_2a}, after unwrapping.}
\end{figure}

In the second example, the orthocomplex exponential is again constant,
but this time `by parts', with two different constant values taken
over the time period considered, as expressed in \eqref{eq:example2_phase}.
The initial frequency of the exponential is 12.5 times that of the highest
frequency of the baseband complex signal, and it then doubles at time
$0.25$ before returning to its original value at time $0.75$.
\figurename~\ref{modulation_2a} shows the modulated signal in this case.
The recovered envelope $\hat{\rho}_2(t)$ is exactly as in the previous example,
the frequency changes in the complex
exponential makes no difference to the recovery of the envelope.
The polar Cayley-Dickson phase $\hat{\phi}_2(t)$ is plotted in \figurename~\ref{modulation_2b} after
unwrapping.

\begin{figure}[tp!]
\includegraphics[width=\columnwidth]{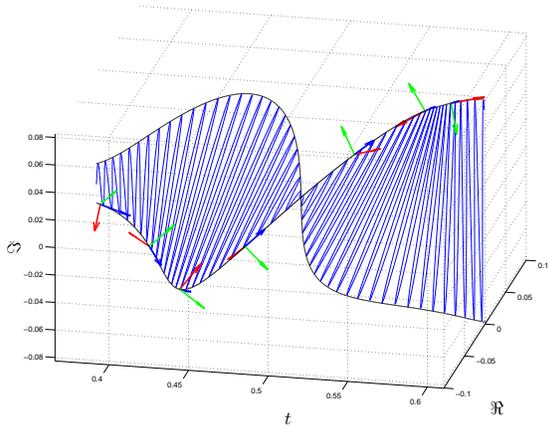}
\caption{\label{modulation_2c}Original signal $z_2(t)$ (blue) and complex
  envelope $\rho_2(t)$(black) of the orthocomplex
  modulation \ref{modulation_2a} for the time period: $t=0.3$ to
  $t=0.6$. For different times, the tangent vectors (blue), normal
  (red) and binormal (green) are displayed. The blue vectors represent
$\hat{\rho}_2'(t)$ for some given $t$, the tangent to the complex
envelope. Green vectors represent $\hat{n}_2(t)$ (given in
\eqref{eq:hatn}), \emph{i.e.} the normal to the instantaneous
osculating plane of the signal $z_2(t)$.}
\end{figure}

\begin{figure}[tp!]
\includegraphics[width=\columnwidth]{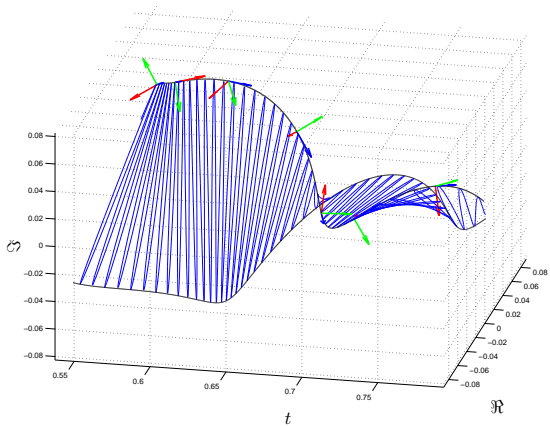}
\caption{\label{modulation_2d}Original signal $z_2(t)$ (blue) and complex
  envelope $\rho_2(t)$(black) of the orthocomplex
  modulation \ref{modulation_2a} for the time period: $t=0.55$ to
  $t=0.8$. For different times, the tangent vectors (blue), normal
  (red) and binormal (green) are displayed. The blue vectors represent
$\hat{\rho}_2'(t)$ for some given $t$, the tangent to the complex
envelope. Green vectors represent $\hat{n}_2(t)$ (given in
\eqref{eq:hatn}), \emph{i.e.} the normal to the instantaneous
osculating plane of the signal $z_2(t)$.}
\end{figure}

For this second example, we also present two zoomed views,
\figurename s~\ref{modulation_2c} and \ref{modulation_2d},
taken for two different time intervals. On
these two figures, in addition to the original signal $z_2(t)$ and the
envelope $\hat{\rho}_2(t)$, some {\em local geometric} information is
displayed. At several times, the vectors $\hat{\rho}_2^\prime(t)$ (blue)
and $\hat{n}_2(t)$ (red) are displayed. One can identify that
$\hat{\rho}^\prime_2(t)$ is the vector tangent to the envelope
$\hat{\rho}_2(t)$  of the signal $z_2(t)$, while the vector
$\hat{n}_2(t)$ is normal to the osculating plane. This local
information gives access to the plane in which the signal is
`oscillating' over time. Note that this \emph{geometric}
information has been obtained from attributes derived from the
hypercomplex representation of the original signal $z_2(t)$, {\em
  i.e.} through simple spectral considerations. Such geometric
quantities are classically obtained from 3D curves using the
Frenet-Serret formulae \cite[Chap. I]{Docarmo1976} for example.

\begin{figure}[ht]
\includegraphics[width=\columnwidth]{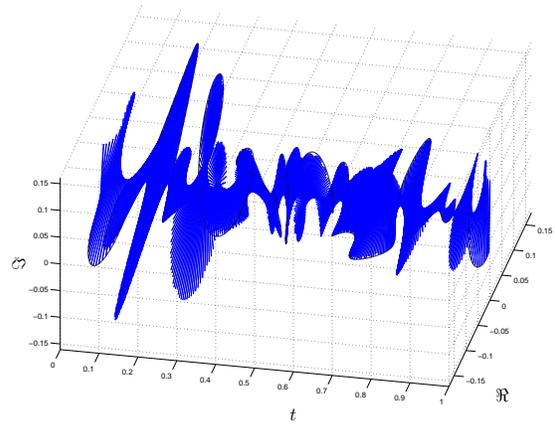}
\caption{\label{modulation_3a}Band-limited signal $z_3(t)$ (blue)  envelope $\hat{\rho}_3(t)$ (black) of the hypercomplex
  representation $\hat{z}_3(t)$ associated to $z_3(t)$.}
\end{figure}

\begin{figure}[ht!]
\includegraphics[width=\columnwidth,height=5cm]{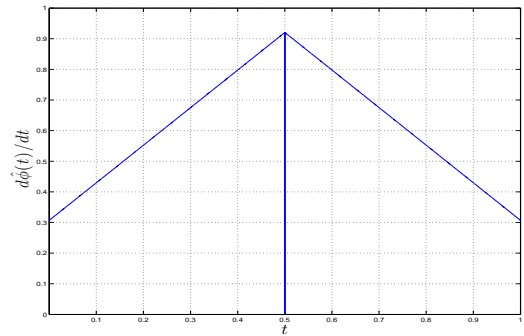}
\caption{\label{modulation_3b}Instantaneous frequency
  $\hat{\phi}_3'(t)$ obtained from the hypercomplex representation
  $\hat{z}_3(t)$ of the complex signal $z_3(t)$.}
\end{figure}

In the third example, the orthocomplex exponential has a much higher frequency
than the baseband signal (which is the same signal as in the first example).
The initial frequency of the exponential is 25 times that of the highest
frequency of the baseband complex signal, and it sweeps linearly to double
the initial frequency, and then back down to the initial value. The
original signal $z_3(t)$ together with the complex envelope $\hat{\rho}_3(t)$ of the
hypercomplex representation $\hat{z}_3(t)$ are displayed in \figurename~\ref{modulation_3a}.
In \figurename~\ref{modulation_3b}, we present
the instantaneous frequency $\hat{\phi}_3'(t)$ obtained from
$\hat{z}_3(t)$. The linearity of the {\em chirp} frequency behaviour
is recovered (singularity at time $0.5$ is due to the behaviour of the
phase at this point) as expected.  

Note that the notion of {\em instantaneous frequency} that we have
introduced here for complex signals is indeed very similar to the
notion of {\em angular velocity} that is well known in physics. The
analogy is possible thanks to the fact that the knowledge of
$\hat{n}(t)$ gives access to the time-varying orientation of the 2D
plane in which the signal is oscillating. 

\begin{figure}[ht!]
\includegraphics[width=\columnwidth]{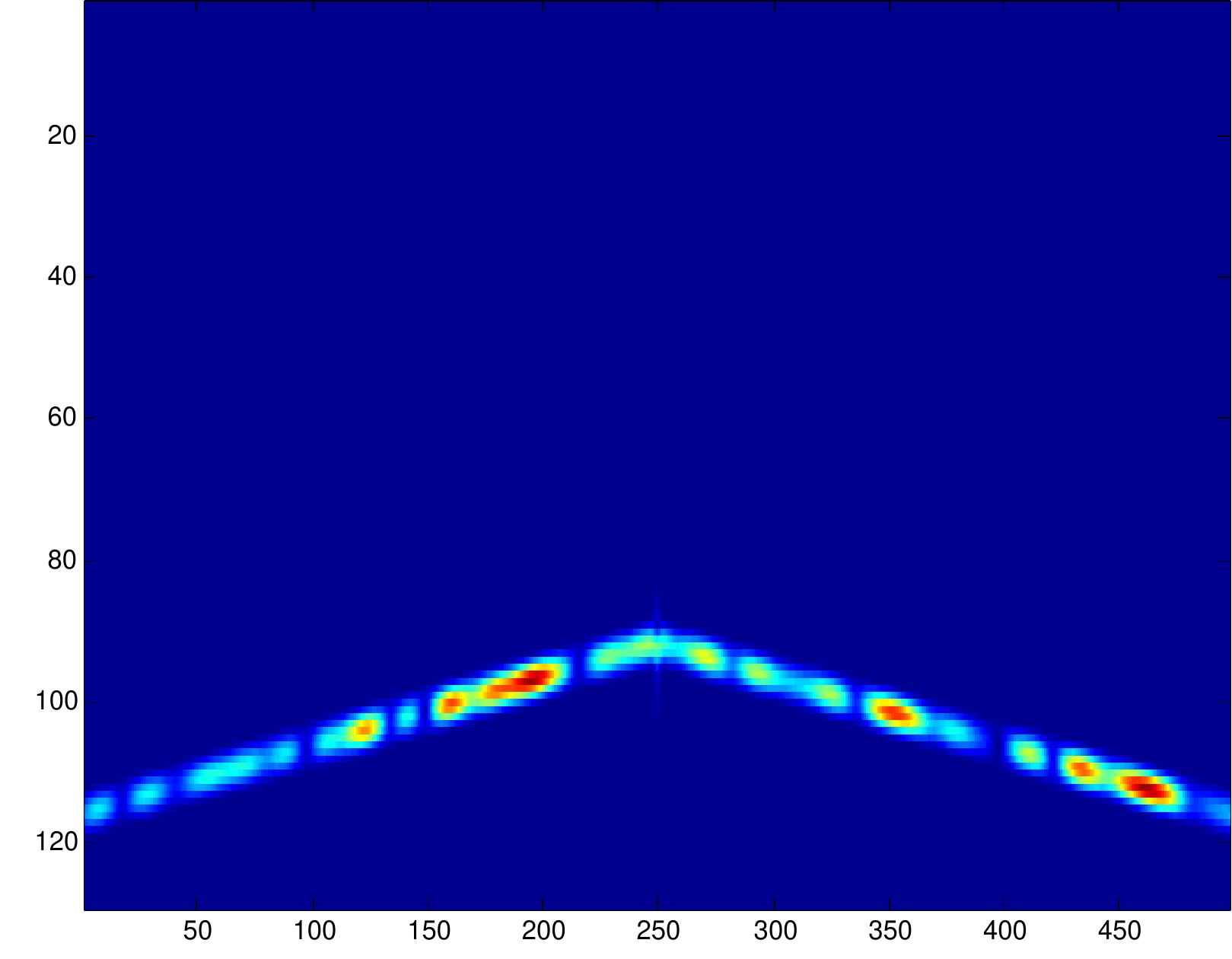}
\caption{\label{timefreq}Time-frequency depiction of the modulated
                         signal $z_3(t)$ in \figurename~\ref{modulation_3a}.
                         The figure shows the modulus of a short-time quaternion
                         Fourier transform with window length 128.}
\end{figure}

Finally, a short-time quaternion Fourier transform (STQFT) of the modulated
signal $z_3(t)$ is shown in \figurename~\ref{timefreq},
showing the frequency sweep of the orthocomplex exponential.
The modulus of the STQFT is real valued and it shows the frequency
content evolution over time for $z_3(t)$. The values for which the
STQFT modulus vanishes represent the moment when the complex envelope
vanishes. 

Notice that we have not illustrated the case where $B_q(t)$ is complex.
This is because it is the \emph{modulus} of $B_q(t)$ which is the phase
in the orthocomplex exponential, and a complex $B_q(t)$ would need to
have a modulus representing phase. It is not clear what the argument
would represent.
A constant frequency is represented by a linearly increasing real phase,
that is $B_q(t)$ would be a real-valued ramp\footnote{In fact it could be
an imaginary-valued ramp with minimal change to the result.}.

As can be seen from \eqref{expBqj}, if $B_q(t)$ is real, that is $B_q(t) = b_o(t)$,
we have the following form for the modulated signal, where $A_q(t) = a_0(t) + a_1(t)\i$:
\begin{equation}
\begin{aligned}
q(t) &= A_q(t)\left[\cos|b_0(t)| + \sin|B_0(t)|\j\right]\\
     &= A_q(t)\cos|b_0(t)| + A_q(t)\sin|B_0(t)|\j
\end{aligned}
\end{equation}
This quaternion signal contains $A_q(t)$ modulated by a cosine in the first
Cayley-Dickson component (the first two components of the quaternion), and
$A_q(t)$ modulated by a sine in the second Cayley-Dickson component (the
third and fourth components of the quaternion).
The orthonormal formulation using quaternions achieves the following:
\begin{itemize}
\item The `carrier' is modeled by a complex exponential, rather than a
      trigonometric function, giving the usual advantages of the
      notation.
\item Because the quadrature component is orthogonal to the in-phase
      component, the two are separated in the signal and in the Fourier
      domain.
\end{itemize}
These advantages disappear if we replace $\j$ with $\i$ and revert to
complex algebra, because then the in-phase and quadrature components
become mixed in a single complex signal, rather than being kept separate
in the components of a quaternion signal.
In complex algebra we can represent a complex exponential modulated
by a real signal, or a real `carrier', represented by a trigonometric
function, modulated by a complex signal. We cannot represent a complex
signal modulating a complex exponential.

\section{Conclusions}
\label{Conclusionsection}
We have shown in this paper how the classical concepts of
instantaneous amplitude and phase may be extended to the
case of complex signals using a quaternion hypercomplex
representation of a complex signal.
We have also shown that the classical case of amplitude
modulation can be extended to the complex case provided
that we maintain an orthogonal separation between the
complex plane of the modulating signal, and the complex
plane of the complex exponential `carrier',
and maintain a separation in frequency between the
modulating signal and the `carrier' exactly as in the
classical case.
The orthogonal separation required is easily realised using
a quaternion representation, and we can then compute a
hypercomplex representation of the complex signal using
a suitable quaternion Fourier transform.
This hypercomplex representation of the signal preserves
symmetries that would be lost using a complex Fourier
transform.
It is analogous to the classical analytic signal.

We have shown that the quaternion polar representation permits us to recover
from the hypercomplex representation, both the instantaneous
amplitude and the instantaneous phase. From the instantaneous phase,
it is also possible to recover the instantaneous frequency as well as
the osculating plane that gives geometric information about the
direction in which the signal is oscillating.  

As quaternions are isomorphic to the geometric Clifford algebra $Cl(0,2)$ of
$\R^2$, it is evident from our results that there are possibilities
to extend the work presented here to the case of vector-valued
signals where the samples are $N$-dimensional,
by using Clifford Fourier transforms in $2^N$ dimensions.

\end{document}